\documentclass{article}
\usepackage[hidelinks]{hyperref}
\usepackage[type={CC},modifier={by-nc-nd},version={4.0}]{doclicense}
\usepackage{fontspec}
\usepackage[scr=rsfs]{mathalpha}

\usepackage[margin=1in]{geometry}

\newenvironment{subproof}[1][\proofname]{%
  \begin{proof}[#1]%
}{%
  \end{proof}%
}

\usepackage{amsmath}
\usepackage{amsthm}
\usepackage{amssymb}

\PassOptionsToPackage{partial=upright}{unicode-math}
\usepackage{fontsetup}
\usepackage{tikz-cd}
\usepackage[style=alphabetic]{biblatex}
\usepackage{braket}
\usepackage{physics}

\title{Triple Products of Eigenfunctions and Spectral Geometry\footnote{Dedicated to Autumn.}}
\date{2026\\September}
\author{Joe Schaefer\\President, SunStar Systems\\ \href{mailto://Joe\%20Schaefer,\%20Ph.D.\%20<joe@sunstarsys.com>}{joe@sunstarsys.com}\footnote{with assistace from Grok.}}

\theoremstyle{plain}
\newtheorem{theorem}{Theorem}
\newtheorem{corollary}[theorem]{Corollary}
\newtheorem{lemma}[theorem]{Lemma}

\theoremstyle{definition}

\newtheorem{remark}[theorem]{Remark}

\begin{document}
\maketitle

\begin{abstract}
Using elementary techniques from Geometric Analysis, Partial Differential Equations, and Abelian $C^*$ Algebras, we uncover a novel, yet familiar, global geometric discriminant -- namely the indexed set of integrals of triple products of eigenfunctions of the Laplace-Beltrami operator, to precisely characterize which isospectral closed Riemannian manifolds are isometric.

\end{abstract}

\tableofcontents

\section{Introduction}
\label{sec:intro}

For a closed, connected Riemannian manifold $(M,g)$, characterizing its class of non-isometric, isospectral manifolds is a type of Inverse Problem \cite{DH11} in Spectral Geometry. Naïvely one might speculate that this class would always be empty. However, the academic literature is rich with decades-old constructions of specific pairings of counterexamples: beginning in 1964 with John Milnor's 16-dimensional pair of non-isometric, isospectral flat tori \cite{JM64}, and continuing \cite{CS92} towards the generic dimensional characterization of flat tori in Alexander Schiemann's 1993 doctoral thesis \cite{AS94} -- replete with a computer aided search for the critical $\dim = 3$ case.  A modern survey of the full flat tori history appears in \cite{NRR22}.

Along the way were insightful offshoots into more sophisticated, non-Euclidean symmetric covering spaces; constructing such isospectral, non-isometric "duets" involving nontrivial curvature tensors (and their spectrum-determined Euler characteristics in dimension 2 \cite{MS67}.) A prime example of this effort was Toshikazu Sunada's 1985 \cite{TS85} invention of a general-purpose covering space framework, which he then deployed in the same work to construct hyperbolic duets in dimensions 2 and 3.

For inhomogeneous Riemannian metrics, Carolyn Gordon discovered duets that are not even locally isometric \cite{CG93}.

Work continues in many related areas \cite{DH11}, such as determining topological characteristics of the class of isospectral, non-isometric manifolds in general (empty \cite{ST80}, finite \cite{AS94}, rigid \cite{GK80}, and compact \cite{GZ97}) as a subset of different moduli spaces of Riemannian metrics.

What we offer in this article is a new perspective on a familiar tool: indexed Fourier coefficients of pairwise products of eigenfunctions as a discrete "algebraic/topological discriminant" to complement the existing, discrete "analytic invariant" -- the non-negative spectrum of the Laplace-Beltrami operator (herein referred to as the (non-negative) \textit{Laplacian}) on $\mathscr H = L^2(M,g)$.  Combined, we observe the pair provides a "discrete global geometric representation" of the isometry classes of isospectral, closed Riemannian manifolds.

If given two closed Riemannian manifolds, and the belief that they were isometric, how would one go about constructing the isometry? This paper is a global roadmap towards addressing that question in fine detail.

This is not a research monograph, but an introduction to the subject matter for advanced math students familiar with Geometric Analysis, Representation Theory, and $C^*$ algebras.

\subsection{Results}
\label{subsec:res}

\begin{theorem}
\label{thm}

Let $(M,g)$ be a closed, connected Riemannian manifold and $\{e^i\}_{i=0}^\infty$ an orthonormal eigenbasis of the non-negative Laplacian $\Delta_M$ on $L^2(M,g)$, ordered by non-decreasing eigenvalues. Define
\begin{equation}
M^{i,j,k} := \int_M e^i e^j \bar{e^k} \sqrt{g} dx = \bra{e^i e^j}\ket{e^k}.
\end{equation}

Another closed Riemannian manifold $(N,h)$ is isometric to $(M,g)$ if and only if it is isospectral to $(M,g)$ and admits an orthonormal eigenbasis of $\Delta_N$ with the same eigenvalues and the same triple products $M^{i,j,k}$.
\end{theorem}

It is important to recognize $M^{i,j,k}$ is not basis-invariant: there is a natural unitary change-of-basis action on it discussed in detail after the proof of Theorem~\ref{thm}. The discussion incorporates certain sets of basis-invariant singular values that one might form a general conjecture around, which claims that the ordered set of singular values completely characterizes the set of isospectral manifolds.

Regardless of the sufficiency half of the general conjecture, necessity is always the case. Which means these collections of singular values defined by $M^{i,j,k}$ and associated to every eigenspace triple are a new set of Riemannian \textbf{invariants}.

The hard work ahead for future research is in locating such basis pairs, or in determining that such pairs cannot exist at all, just by examining the properties of the $M^{i,j,k}$ in evidence. But this paper puts that target front and center: we seek to reduce the analytic geometry questions of Spectral Theory to computationally tractable linear algebra questions about products of eigenfunctions.

\textit{Symmetry} plays an important rôle in computationally tractable cases \cite{TF17} \cite{LS18} \cite{PS94}, which is aptly illustrated in our flat tori Example~\ref{sec:ex} below. However, the strength of our approach is perhaps best made apparent in the case of manifolds with the fewest number of Riemannian symmetries, which is the generic case. In this instance, we offer the following

\begin{corollary}
\label{coro-1}
(Diagonal Litmus Test) Given a pair of eigenvalue preserving orthonormal bases as described in the hypothesis of the Theorem, the manifolds are isometric if for every choice of $i,j,k$, the product $M^{i,\bar i,k}\bar M^{j,\bar j,k}$ agrees in both bases; and if the vector space spanned by $\set{|e^i|^2}$ is dense in $\mathscr H$. Here $\bar j$ represents the eigenfunction $\bar e^j$ in the triple-product integral computations.

Furthermore, if we define $\mathscr V$ as the Hilbert space generated by $\set{|e^i|^2}$, $\mathscr V = \mathscr H$ if and only if the adjoint map

\begin{equation}
[M^{i,\bar i, k}]^*:\mathscr H \rightarrow \mathscr V
\end{equation}

is injective.
\end{corollary}

\begin{corollary}
\label{coro-2}
Generically, isospectral manifolds are isometric if and only if the products as defined in Corollary~\ref{coro-1} agree as real values.
\end{corollary}

The motivation for the study of $\set{M^{i,j,k}}$ is loosely derived from the study of the rôle of the bilinear multiplication operator $Y:V\otimes V\rightarrow V((z))$ in the definition of a Vertex Operator Algebra \cite{FBZ04} associated with a Chiral Conformal Field Theory. Here $V$ is the Vector Space of States and $V((z))$ is the space of formal Laurent series in $z$ with coefficients in $V$. Since $V$ often comes equipped as a Hilbert Space with a traditional Fourier series orthonormal basis, indexing $Y$ using the Fourier basis elements of $V$ is only slightly more involved than the $M^{i,j,k}$ case studied here, but quite similar in spirit. However a detailed comparison is out of scope for this article.

If we consider the map
\begin{equation}
(M, g, \set{e^i}) \mapsto \set{\lambda_i, M^{i,j,k}}\ ,
\end{equation}
this paper establishes the injectivity of this map for closed Riemannian manifolds (up to Riemannian isometry in its domain). Further results which apply these techniques are just getting started. Adve \cite{AA25} proves a converse of a different kind: a discrete solution of the $\mathrm{PSL}_2(\mathbb{R})$ crossing equations arises from a compact hyperbolic 2-orbifold. Crossing is associativity in a $G$-isotypic basis (Remark 5). It is not an input to Theorem 1.

Some imagery may be helpful here. If we fix $(M,g)$ and look at the orbits of $\set{M^{i,j,k}}$ under spectrum-preserving change-of-basis unitary transformations on $\set{e^i}$, we see that the orbits of different isospectral $(M,g)$ pairs partition the image of this map along isometry classes.

Finally, we prove that the generic Riemannian metric case is completely characterized by the study of the "diagonal" $\set{\lambda_i, M^{i,\bar i,k}}$.

These results were first demonstrated during a similarly titled talk by the author at \textbf{MSRI} in 1997, but they appear here in published form for the first time.

\section{Preliminaries}
\label{sec:pre}

Now with $M,g,e^i,M^{i,j,k}$ as in the hypotheses of Theorem~\ref{thm}, for $f \in C^\infty(M)$ and $i \geq 0$ note that the \textit{Fourier coefficients}

\begin{equation}
\label{eqn:Fourier}
\begin{aligned}
\hat{f}(i) &:= \int_M f(x)\bar{e^i}(x)\sqrt{g(x)}dx \\
\implies \\
f(x)        &= \sum_{i=0}^{\infty}\hat{f}(i)e^i(x)\,,
\end{aligned}
\end{equation}
since $f$ is uniquely representable as its rapidly converging \textit{Fourier Series} ($\Delta_M$-specific Sobolev Embeddings \cite{MT13} \cite{RS75}, together with Weyl's Asymptotic Law \cite{HW11}, imply the terms in the sum are $o(i^{-n})$ uniformly in $x$ \cite{LH68}, $\forall n\in\mathbb{N}$.) Then we see that for $f_1, f_2 \in C^\infty(M)$, the Fourier coefficients of the pointwise product $f_1 f_2 \in C^\infty(M)$ are

\begin{equation}
\begin{aligned}
\widehat{f_1 f_2}(k) &= \sum_{i,j}^\infty\hat{f_1}(i)\hat{f_2}(j)M^{i,j,k} \\
\implies \\
f_1f_2(x) &= \sum_{i,j,k}\hat{f_1}(i)\hat{f_2}(j)M^{i,j,k}e^k(x) \\
f_1 &= f^p_2,\, p > 2 \implies \\
\sum_{k}\hat{f_1}(k)e^k(x) &= \sum_{i_1,i_2,\ldots,i_{2p-1}}\hat{f_2}(i_1)\hat{f_2}(i_2)\hat{f_2}(i_4)\hat {f_2}(i_6)\ldots\hat{f_2}(i_{2p-2})M^{i_1,i_2,i_3}M^{i_3,i_4,i_5}\ldots M^{i_{2p-3},i_{2p-2},i_{2p-1}}e^{i_{2p-1}}(x)
\end{aligned}
\end{equation}
and so, critically, any multivariate polynomial $\wp \in \mathbb{C}[z_1,\ldots,z_l]$ (on smooth functions) commutes with any spectrum-preserving $\Delta$-eigenfunction orthonormal basis map $\vec{F}$ that preserves $\set{M^{i,j,k}}$:

\begin{equation}
\label{fig}
\begin{tikzcd}
  C^\infty(M,\space\mathbb{C}^l)\arrow[swap]{d}{\underbrace{\vec{F}\oplus\dots\oplus \vec{F}}_{l\,\text{times}}} \arrow{r}{\wp} & C^\infty(M) \arrow{d}{\vec{F}}\\%
  C^\infty(N,\space\mathbb{C}^l) \arrow{r}{\wp} & C^\infty(N)
\end{tikzcd}
\end{equation}

Moreover if $A\subset M$ is Borel-measurable, then the results above hold pointwise for the characteristic function of $A$ everywhere except along the boundary of $A$: if $f = f^2$ and $A:=\set{x\in M|f(x)=1}$,

\begin{equation}
\sum_{i}\hat{f}(i)e^i(x) = \sum_{i,j,k}\hat{f}(i)\hat{f}(j)M^{i,j,k}e^k(x) = \begin{cases}
1 & x \in \mathring{A} \\
0 & x \in \mathring{A^\complement}\end{cases}\,,
\end{equation}
and by uniqueness, we have the following identity

\begin{equation}
\begin{aligned}
\hat{f}(k) &= \sum_{i,j}\hat{f}(i)\hat{f}(j)M^{i,j,k}\,,\, \forall k\geq 0 \\
\iff f&=f^2 \, a.e.
\end{aligned}
\end{equation}

This implies any such basis map as above carries characteristic functions (as members of $L^2(M,g)\subset L^1(M,g)$) to characteristic functions in a measure-preserving fashion.

The point of these computations is to emphasize the fact that $\set{M^{i,j,k}}$ characterizes the Harmonic Analysis of the pointwise multiplication operator on $C^\infty(M)$, which is a dense subalgebra of the Abelian $C^*$ algebra $C(M)$, by the Stone-Weierstrass theorem.

For the rapid convergence of these above sums involving $M^{i,j,k}$, note that products of eigenfunctions are smooth, so these Fourier coefficients decay as above (in each index). For more details, see Emmett Wyman's work in 2022 with these coefficients as it relates to the triangle inequality on the eigenvalues \cite{EW22}.

Note: we may always assume

\begin{equation}
\begin{aligned}
e^0 &= M^{0,0,0} = 1/\sqrt{vol(M)} \\
\implies \\
M^{0,j,k} &= M^{j,0,k} = \delta_{j-k}\space/\sqrt{vol(M)}\,,
\end{aligned}
\end{equation}
where $\delta_i$ is the Kronecker delta. Since $vol(M)$ is a spectral invariant \cite{HW11}, this information is already available from isospectrality considerations.

\section{Proof of Theorem~\ref{thm}}
\label{sec:proof}

\begin{proof}
For necessity, let $F:(N,h)\rightarrow (M,g)$ be an isometry between closed Riemannian manifolds, and let the target orthonormal basis of eigenfunctions on $L^2(N,h)$ be the pull-back via $F$ of the orthonormal basis $\set{e^i}$ on $(M,g)$ above. Since

\begin{equation}
\begin{aligned}
M^{i,j,k} &= \int_M e^i e^j \bar{e^k}\sqrt{g}dy \\
          &= \int_N e^i(F(x)) e^j(F(x))\bar{e^k}(F(x))\sqrt{h}dx\,,
\end{aligned}
\end{equation}
we are done with the necessity argument because $\Delta_N(f\circ F) = (\Delta_M f) \circ F,\ \ \forall f\in C^\infty(M)$.

For sufficiency, we now consider the linear, bijective orthonormal eigenfunction basis map $\vec{F}$ from $C^\infty(M)$ to $C^\infty(N)$ and note that from the calculations in Section~\ref{sec:pre} above, $\vec{F}$ preserves pointwise products for smooth functions (and preserves characteristic functions when extended to $L^2(M,g)$) by the premise that $\set{M^{i,j,k}}$ is invariant under this map.

\begin{lemma}
\label{lma}
$\vec{F}: C^\infty(M)\rightarrow C^\infty(N)$ preserves the uniform norm.
\end{lemma}

\begin{subproof}[Proof of Lemma~\ref{lma}]
Continuity of $\vec F$ is a consequence of the Sobolev embedding theorem and the fact that $\vec F$ is an isometry on every Sobolev space between $M$ and $N$. Since $\vec F$ is bijective, we just need to establish $$\lVert\vec F(f)\rVert_\infty = \lVert f\rVert_\infty.$$ From the Preliminaries section above, we can compute
$$
  \lVert\vec F(f)\rVert_{2p}^{2p}=\lVert\vec F(f^p)\rVert_2^2 = \lVert f^p \rVert_2^2 = \lVert f \rVert_{2p}^{2p},
$$
and so taking $p\to\infty$ proves the Lemma. However, we will take a slightly different tack in the proof provided below.

Let $\set{a_i}$ be a smooth partition of unity on $M$.

\begin{equation}
\begin{aligned}
1 &= \sum_i a_i(x) \\
       &= \sum_{i,j} \hat{a_i}(j)e^j(x) \\
	   &= \sum_j e^j(x)\sum_i \hat{a_i}(j)
\end{aligned}
\end{equation}

Thus $\sum_i\hat{a_i}(j) = \delta_j\sqrt{vol(M)}$ (Kronecker delta).

By the dominated convergence theorem,
\begin{equation}
\lim_{p\rightarrow\infty}
\sum_j\hat{a^p_j}(k) = \int_{\dot{\bigcup}_j\set{a_j=1}}\bar{e^k}(x)\sqrt{g}dx\,,
\end{equation}
which the $k$-th Fourier coefficient of a characteristic function (assumed of positive measure) on each disjoint subset $\set{x\in M | a_j(x) = 1}$. This means the Lemma is proven for each $a_j$, since the limiting characteristic function of a set with positive measure is preserved, and hence has uniform norm 1, as do all $a_j^p,\, \vec{F}(a_j^p)=\vec{F}(a_j)^p,\, p\in\mathbb{N}$, by Diagram~\ref{fig}.

Without loss of generality, we may apply the special case result shown for  the smooth partition of unity $\lbrace|f|/\lVert f \rVert_\infty, 1 - |f|/\lVert f\rVert_\infty\rbrace$, where $ \set{x\in M|\space|f(x)| = \lVert f \rVert_\infty}$ has positive measure, and the Lemma is proven in full.

\end{subproof}

Since $\set {\bar e^i}$ is also a Fourier basis for $L^2(M,g)$, it is clear from Equation~\ref{eqn:Fourier} that $\vec F(\bar f) = \bar{\vec F}(f)$. This means that on a dense set of $C(M)$ (and $C(N)$), we have established $\vec{F}$ as an isomorphism of Abelian $C^*$ algebras, and thus can be extended to an isomorphism of $C(M)$ and $C(N)$ in the same category.

Now we apply the Gelfand-Naimark Representation Theorem (in contravariant functor form) for unital Abelian $C^*$ algebras \cite{JC19} to represent this isomorphism $\vec{F}$ by a homeomorphism $F$ between $N$ and $M$. Since $\vec{F}$ is bijective on smooth functions, $F$ too must be smooth.

As this now diffeomorphism preserves eigenvalues and eigenfunctions (by hypothesis on $\vec{F}(f) = f \circ F$), it must preserve the Laplacian on smooth functions. Hence $F$ also must preserve the principal symbols of these same elliptic operators \cite{MT13}. The principal symbols of the Laplacian are simply another means of expressing the Riemannian metric on the manifolds in question.

This completes the proof of the Theorem.

\end{proof}

\subsection{Discussion of Corollaries}
\label{subsec:coros}

With $\set{M_0^{i,j,k}}$ and $\set{M_1^{i,j,k}}$ representing the two triple-product sets for the bases $\set{e_0^i}$ and $\set{e_1^i}$, let $z_i \in U_1$ be the $U_1^\infty$ action on such an orthonormal basis $\set{e_1^i}$. Thus, we will choose $z_i$ so that $\set{z_ie_1^i}$ yields $\set{M_0^{i,j,k}} = \set{z_i z_j \bar z_kM_1^{i,j,k}}$.

Why is this the case? In general, the symmetry group acting on the space of possible orthonormal bases of eigenfunctions is the space of Unitary Operators $U: \mathscr H\rightarrow\mathscr H$ that commute with projections $P_{\mathscr V_\lambda}$ onto the finite-dimensional eigenspaces $\mathscr V_{\lambda}$ associated with each individual eigenvalue $\lambda$ of the Laplacian. Therefore

\begin{equation}
\begin{aligned}
P_{\mathscr V_{\lambda}}U(e^i) = UP_{\mathscr V_{\lambda}}(e^i),\ \therefore U(e^i) &= \sum_{\lambda_i = \lambda_j}u_{ij}e^j \implies \\
M_U^{i,j,k} := \int_M U(e^i)U(e^j)\bar U(\bar e^k)\sqrt g dx &= \sum_{\lambda_r = \lambda_i,\lambda_s=\lambda_j,\lambda_t=\lambda_k} u_{ir}u_{js}\bar u_{tk} M^{r,s,t}
\end{aligned}
\end{equation}
is the image of $M^{i,j,k}$ under $U$’s basis action $e^i \mapsto U(e^i)$.

Now under the conditions of Corollary~\ref{coro-2}, each of the $\mathscr V_\lambda$ are one dimensional vector spaces over $\mathbb{C}$, but that also means they are one dimensional vector spaces over $\mathbb{R}$, and so the full multiplicative symmetry group is $O(1,\mathbb{R})^\infty=\mathbb{Z}_2^\infty$.

More generally, the Corollary~\ref{coro-2}'s associated prerequisite "regarding agreement in product values" would simply become "preservation of the ordered set of singular values (counted with multiplicity) of the linear map $M_{\lambda_i\lambda_j\to\lambda_k}:\mathscr V_{\lambda_i} \otimes \mathscr V_{\lambda_j}\to \mathscr V_{\lambda_k}$ defined by $\set{M^{i,j,k}}$ on the basis elements, with the Hilbert-Schmidt norm on the tensor-product domain. By definition, these singular values are invariant under direct sums of unitary transformations on the $\mathscr V_\lambda$.

In the multiplicity-1 spectrum case, the complete set of singular values is simply the set of absolute values of $M^{i,j,k}$ which, we still conjecture, completely characterizes the isometry classes of such isospectral manifolds. See Equation~\ref{assoc} for the key relationship between this conjecture and Corollary~\ref{coro-2}. What's missing is the sufficiency argument that if the absolute values agree, the manifolds are isometric; which requires an argument to eliminate possible sign change cancellations between bases in the LHS summands of Equation~\ref{assoc}.

We are significantly less confident that the general conjecture holds true (outside the multiplicity-1 spectrum case), since it may be possible to produce a counterexample (of sufficiency) via explicit Sunada construction.

If the index notation is obfuscating the situation, perhaps this basis-independent description will help.  Take $v_{\lambda} \in \mathscr V_{\lambda}$ and consider the expression

$$
P_{\mathscr V_\gamma}(v_\alpha v_\beta).
$$

Recall that each $\mathscr V_\lambda$ is a finite dimensional complexified Euclidean space. All $M^{i,j,k}$ does is provide basis coordinates of this basis-independent expression. Since $1 = \oplus_\lambda P_{\mathscr V_\lambda}$, what the Theorem says is that the above expression is \textbf{identical} between manifolds if and only if the manifolds are isometric; which should come as a shock to literally no one.

The rôle of the $M^{i,j,k}$ is purely relative (comparative). They give us a way of comparing these tensor product operators between two different manifolds. Why do we care? Because we want to use these coordinates to justify the promotion of the associated basis map to a Riemannian isometry between the manifolds.

The sufficiency half of these conjectures are largely combinatorics problems involving reconstructing these expressions purely from their singular value decomposition.

However these basis invariants may prove useful in deciphering complex cases to prove two manifolds are \textit{not isometric}, by showing that their singular values are not identical between the two bases in question.

\subsection{Aside}

The Representation Theory of a Compact Lie Group $G$ takes the explicit Laplacian out of the equation and studies \textit{irreducible} representations $\Phi^{(\alpha)} \in Hom(G,U(\mathscr V_\alpha))$ towards the Hilbert space decomposition of $L^2(G,dg)$ (here $dg$ is the normalized Haar Probability Measure on $G$) as addressed in the Peter-Weyl Theorem, and honors their interplay in the expression above (via irrep decomposition of tensor products of irreps) as the essential artifacts of Lie Theory \cite{AK01}. Compatible Riemannian geometries are generated by convenient choices of quadratic Casimir elements that lie in the center of the universal enveloping algebra, which are of less significance than the irrep matrix coefficient decomposition of $L^2(G,dg)$ itself. \textbf{Their spectral decomposition} is a (less convenient) reassembly of those orthonormal basis elements, since the associated Casimir element is constant on each irrep's matrix coefficients.

Wigner's $3j$ symbols for $SU(2)$ are a prime example for further study -- let us apply our general construction to it as a group manifold.  Every finite-dimensional irreducible unitary representation of $SU(2)$ is labelled by a non-negative half-integer $j = 0, \tfrac12, 1, \tfrac32, \dots$. We write $\mathscr  V_j$ for the $(2j+1)$-dimensional space on which this representation acts. A standard orthonormal basis of $\mathscr  V_j$ is the magnetic basis $\ket{j\ m}$ with $m = -j,-j+1,\dots,j$.

The equation

$$
\begin{pmatrix}
j_1 & j_2 & j_3 \\
m_1 & m_2 & m_3
\end{pmatrix}
:=(-1)^{j_1-j_2-m_3}\frac{1}{\sqrt{2j_3+1}}C^{j_3\ -m_3}_{j_1\ m_1,j_2\ m_2}
$$
expresses the $3j$ symbol definition in terms of Clebsch-Gordan Coefficients $C^{j\ m}_{j_1\ m_1,j_2\ m_2}:=\bra{j_1\ m_1,j_2\ m_2}\ket{j\ m}$ that decompose the tensor product representation $\mathscr V _{j_1}\otimes\mathscr V _{j_2}$ into irreducible components, which have closed form expressions such as Racah's formula that underpins modern numerical software libraries \cite{JF16}.

A $3j$ symbol vanishes unless $m_1+m_2+m_3=0$, the triangle inequalities $|j_1-j_2|\le j_3\le j_1+j_2$ hold, and $j_1+j_2+j_3$ is an integer.

From the point of view of the compact group $SU(2)$, the $3j$ symbols are precisely the (properly normalized and phased) intertwining operators that realize the unique (up to scale) invariant subspace of the triple tensor product $\mathscr V_{j_1}\otimes \mathscr V_{j_2}\otimes \mathscr V_{j_3}$ when that product contains the trivial representation. They are therefore the natural "structure constants" for the fusion of three irreducible representations to the singlet.

This is exactly analogous to the rôle played by the triple-product integrals $M^{i,j,k}$ on a Riemannian manifold: they are the structure constants of the pointwise product of eigenfunctions when that product is expanded back in the eigenbasis. On the group $SU(2)$ itself (or on $S^2=SU(2)/U(1)$) those integrals reduce to phased products of the $3j$-symbols.

To wit, let $D^{j}_{m n}(g)$ be the standard Wigner $D$-functions (matrix coefficients of the irrep of spin $j$ in the magnetic orthonormal basis $\ket{j\ m}$). With respect to the normalized Haar probability measure $dg$ one has the exact formula \cite{VK88}:

\begin{equation}
\int_{SU(2)}
D^{j_1}_{m_1 n_1}(g)\,
D^{j_2}_{m_2 n_2}(g)\,
\overline{D^{j_3}_{m_3 n_3}(g)}\,dg
=
\begin{pmatrix}
j_1 & j_2 & j_3 \\
m_1 & m_2 & -m_3
\end{pmatrix}
\begin{pmatrix}
j_1 & j_2 & j_3 \\
n_1 & n_2 & -n_3
\end{pmatrix}
\times
(-1)^{m_3+n_3}.
\end{equation}
(The overall phase convention can be adjusted by the usual Condon–Shortley factors; the essential point is that the integral factors into a product of two real $3j$-symbols.) The Peter–Weyl theorem supplies a complete orthonormal basis of $L^2(SU(2),vol(G)dg)$ given by the renormalized matrix coefficients:

$$e^{j,m,n}(g):=\sqrt{\frac{2j+1}{vol(G)}}\ D^{j}_{mn}(g),$$
where the indices run over

$$j=0,\tfrac12,1,\tfrac32,\dots,\ m,n=-j,-j+1,\dots,j.$$

Of course, under the Cartan-Killing Casimir Laplacian, $\Delta e^{j,m,n} = \lambda_j e^{j,m,n} = j(j+1)e^{j,m,n}$, which has $vol(G)=2\pi\ vol(S^2) = 8\pi^2$ under its notion of (metrized) Haar measure: $\sqrt g dx = 8\pi^2 dg$. Combining the two equations while sticking to the conventions of the Theorem, one obtains the explicit expression:

\begin{equation}
\label{eqn-su2}
\begin{aligned}
M^{(j_1 m_1 n_1),(j_2 m_2 n_2),(j_3 m_3 n_3)}
&=
\sqrt{\frac{(2j_1+1)(2j_2+1)(2j_3+1)}{8\pi^2}}\\
&\ \times
\begin{pmatrix}
j_1 & j_2 & j_3 \\
m_1 & m_2 & -m_3
\end{pmatrix}
\begin{pmatrix}
j_1 & j_2 & j_3 \\
n_1 & n_2 & -n_3
\end{pmatrix}
(-1)^{m_3+n_3}.
\end{aligned}
\end{equation}

Thus, we see that the Riemannian Geometry of $SU(2)$ induced by the Cartan-Killing Casimir can be captured using our augmented spectral data $\set{\lambda_i, M^{i,j,k}}$, from which singular value Riemannian invariants can be explicitly computed from Equation~\ref{eqn-su2}:

$$\sigma(j_1,j_2,j_3)
=
\sqrt{\frac{(2j_1+1)(2j_2+1)(2j_3+1)}{8\pi^2}}\cdot\frac{1}{2j_3+1}
=
\sqrt{\frac{(2j_1+1)(2j_2+1)}{(2j_3+1)\,8\pi^2}},
$$
whenever the triangle inequalities and parity condition are satisfied on $j_1,j_2,j_3$, with multiplicity 1; and all other singular values zero. The computation follows from the observation that $M^{i,j,k}$ decomposes into the tensor product of left/right rank-1 intertwining operators, each with Hilbert-Schmidt norm $1/\sqrt{2j_3+1}$ (by definition), so its non-zero singular value, per definition above, is the Peter-Weyl rescaled product of those Hilbert-Schmidt norms.

For example, $  (j_1,j_2,j_3)=(1,1,0)  $ the 3 non-zero $  3j  $-symbols have absolute value $  1/\sqrt{3}  $, the prefactor is $  \sqrt{3\cdot3\cdot1/8\pi^2}  $, and one recovers$$\sigma(1,1,0)=\sqrt{\frac{9}{8\pi^2}}\cdot 1=\frac{3}{\sqrt{8\pi^2}},$$which matches the direct Hilbert–Schmidt computation.

However, since $SU(2)$ with the Cartan-Killing Metric is isometric to $4$ times the traditional homogeneous metric on $S^3$, and for $S^n, n\leq 6$, every isospectral manifold is isometric \cite{ST80} (a beautiful reduction to constant sectional curvature metrics from the short-time asymptotics of the trace of the heat kernel), this computation of the singular values is largely academic for $SU(2)$.

Since the homogeneous metric case for $n>6$ is presently unsettled, performing the same Clebsch-Gordan analysis on $SO(n+1)$ \cite{BJK78} with the aid of a computer \cite{HR11} \textit{might} lead to new discoveries about the spectral determinism of $S^n$ in terms of its singular values \cite{M06}.

For compact Abelian Lie Groups, these vector spaces $\mathscr V_\alpha$ are all one-dimensional, so their situation is entirely similar to the spectral decomposition of multiplicity-1 Laplacians above.  More on this in the Example below.

Getting back to the Corollary~\ref{coro-1}, we observe that the proof involves establishing this implication:

\begin{equation}
\label{fml}
z_k = M_0^{i,\bar i,k} / M_1^{i,\bar i,k} \,\, \forall i,k\in\mathbb{N},\, ⋺ M_0^{i,\bar i,k} \ne 0 \, \implies
\exists r,s,t \in \mathbb{N}\ ⋺\ \frac{M_0^{i,j,k}}{M_1^{i,j,k}} = \frac{M_0^{r,\bar r,i}M_0^{s,\bar s,j}\bar M_0^{t,\bar t,k}}{M_1^{r,\bar r,i}M_1^{s,\bar s,j}\bar M_1^{t,\bar t,k}}\,.
\end{equation}
We may hope that for any given $k>0$, $M_0^{i,\bar i,k}$ cannot be identically $0$ for all $i$, since it is a generically true condition, but false for specific cases like the flat tori case covered in the Example below. A higher-level way of looking at this condition is to note that such a hope-violating $k$ would have $\bar e^k$ in the kernel of the adjoint map $[M^{i,\bar i,k}]^*$. Furthermore, the Formula~\ref{fml} for $z_k$ requires both $i$-independence, and sufficiency, to establish the basis map $e_0^i \mapsto z_i e_1^i$ preserves $\set{M_0^{i,j,k}}$.

We sketch a proof of Corollary~\ref{coro-1} (sufficiency) below the next set of formulae.

Nevertheless, let us compute some relevant identities so some intrepid future researcher can dig into this generalized conjecture. Here $v\cdot w$ is the Riemannian cometric:

\begin{equation}
\label{eqn-15}
\begin{aligned}
\Delta fg &= f\Delta g + g\Delta f - 2 df \cdot dg \implies \\
M^{i,j,k} &= 2 \frac{\bra{de^i\cdot de^j}\ket{e^k}}{\lambda_i +\lambda_j -\lambda_k}\\
\text{and so the quadratic form} \\
Q_k(f,g) :&= \bra{df\cdot dg}\ket{e^k} = \sum_{i,j}\hat{f}(i)\hat{g}(j)\bra{de^i\cdot de^j}\ket{e^k} \\
&= \frac{1}{2}\sum_{i,j}\hat{f}(i)\hat{g}(j)(\lambda_i + \lambda_j - \lambda_k)M^{i,j,k}.\\
\Gamma(f,g) := df \cdot dg &= \sum_k Q_k(f,g)e^k = -\frac{\Delta fg - f\Delta g - g\Delta f}{2}\\
\text{And so we have the Bochner Formula:}\\
\Gamma_2(e^i,e^j) := \frac{1}{4}\sum_k (\lambda_i + \lambda_j - \lambda_k)^2M^{i,j,k}e^k&=\langle\operatorname{Hess}e^i,\operatorname{Hess}e^j\rangle_{\mathrm{HS}}+\operatorname{Ric}(\nabla e^i,\nabla e^j),\\
\Gamma_\ell(e^i,e^j) :&= \frac{1}{2^\ell}\sum_k (\lambda_i + \lambda_j - \lambda_k)^\ell M^{i,j,k}e^k.
\end{aligned}
\end{equation}

Cauchy–Schwarz on the Hessian, $\lvert\operatorname{Hess}f\rvert_{\mathrm{HS}}^2\ge\frac1n(\Delta f)^2$, turns $\operatorname{Ric}\ge K\,g$ into exactly $\mathrm{BE}(K,n)$ \cite{BE85}. The bound is sharp on space forms. Integrating these inequalities from the Bochner Formula for $\Gamma_2(e^1,e^1)$ yields the Lichnerowicz lower bound for $\lambda_1$ (only $k=0$ survives integration):

$$\lambda_1\ge\frac{n}{n-1}K.$$

The Hessian is not an extra datum. In $\Gamma$-calculus it is
$$2\operatorname{Hess}f(\nabla u,\nabla v)
=\Gamma\bigl(u,\Gamma(f,v)\bigr)+\Gamma\bigl(v,\Gamma(f,u)\bigr)-\Gamma\bigl(f,\Gamma(u,v)\bigr).$$
Each $\Gamma(\,\cdot\,,\,\cdot\,)$ is a function whose eigen-coefficients are known, so the right-hand side is a finite combination of $M$’s and $\lambda$’s. Concretely, write
$$\gamma_{ab}^c:=\frac12(\lambda_a+\lambda_b-\lambda_c)M^{a,b,c},$$
so $\Gamma(e^a,e^b)=\sum_c\gamma_{ab}^c e^c$. Then
$$\Gamma\bigl(e^p,\Gamma(e^i,e^q)\bigr)
=\sum_{c,r}\gamma_{iq}^c\gamma_{pc}^r e^r,$$
and
$$\begin{aligned}
2\operatorname{Hess}e^i(\nabla e^p,\nabla e^q)
&=\sum_r\Bigl(
\sum_c\gamma_{iq}^c\gamma_{pc}^r
+\sum_c\gamma_{ip}^c\gamma_{qc}^r
-\sum_c\gamma_{pq}^c\gamma_{ic}^r
\Bigr)e^r.
\end{aligned}$$
The Hilbert–Schmidt pairing $\lvert\operatorname{Hess}e^i\rvert_{\mathrm{HS}}^2$ is obtained by tracing this $2$-tensor against the cometric. The cometric itself is $\Gamma$:
$$du \cdot dv =\Gamma(u,v).$$
Generically, finitely many eigenfunctions give local coordinates, their $\Gamma$-matrix is $g^{ab}$, and one inverts to $g_{ab}$. Then
$$\operatorname{Ric}(\nabla e^i,\nabla e^j)
=\Gamma_2(e^i,e^j)-\langle\operatorname{Hess}e^i,\operatorname{Hess}e^j\rangle_{\mathrm{HS}}$$
is an explicit (if cumbersome) series in $\{\lambda_a,M^{a,b,c}\}$.

\section{Proof of Corollaries}
\label{sec:proof-of-coros}
\begin{proof}

Now consider the famous associativity relations from Conformal Field Theory:

\begin{equation}
\label{assoc}
\begin{aligned}
e^ie^je^k = \sum_\ell\bra{e^ie^j}\ket{\bar{e^k}e^\ell}e^\ell &= \sum_{\ell,r} M^{i,j,r}\bar M^{\bar k,\ell,r}e^\ell\\
&= \sum_{\ell,r} M^{i,k,r}\bar M^{\bar j,\ell,r}e^\ell\ \therefore\\
i = \bar j, \ell = k\text{ and relabeling } \implies \\
\sum_r M^{i,\bar i,r}\bar M^{j,\bar j,r} &= \sum_r |M^{i,j,r}|^2
\end{aligned}
\end{equation}

Corollary~\ref{coro-1} follows from the fact that $z_k\in U_1$ is well-defined under its presumed hypotheses, and the prior observation that the algebraic, bounded, trilinear operators defined by $z_iz_j\bar z_kM_1^{i,j,k}$ and $M_0^{i,j,k}$ are both associative, and agree with pointwise function multiplication by squares of absolute values of eigenfunctions, the linear span of which is dense in $\mathscr H$.  Establishing $\ker\ [M^{i,\bar i,k}]^* = 0$ is exactly equivalent, where $\mathscr V$ is the closed Hilbert space generated by $\set{|e^i|^2}$, and $[M^{i,\bar i,k}]:\mathscr V\rightarrow \mathscr H$ is the change-of-basis identity map.

So they agree everywhere. This completes the proof of Corollary~\ref{coro-1}.

Corollary~\ref{coro-2} sufficiency follows by noting that the vanishing adjoint kernel condition in Corollary~\ref{coro-1} is generically true. And if for some choice of $i,j,k$, the product $M^{i,\bar i,k} \bar M^{j,\bar j,k}$ disagreed between bases, it would disagree in every pair of bases.

Why? Since generic manifolds can be presumed to also have multiplicity-1 spectra, this reduces the full symmetry group to $U_1^\infty$ where these products are invariants, contradicting Theorem~\ref{thm}. Further reduction to $\mathbb{Z}_2^\infty$ via real-valued bases ensures the products are real-valued. This establishes the necessity of the hypothesis, and completes the proof of Corollary~\ref{coro-2}.

Further, the arguments in the proof of Corollary~\ref{coro-1} are valid even when the manifolds are non-isospectral, so we can represent the basis map as a homeomorphism if and only if the products match in the generic case above. Isospectrality then becomes equivalent to this homeomorphism being promoted to a Riemannian isometry -- due to Sobolev norm equivalence on the corresponding Fourier coefficients (realizing it as a diffeomorphism), and leading to principal symbol preservation on the Laplacians for full isometry as in the proof of the Theorem.

\end{proof}

\begin{remark}[Associativity and Crossing] Suppose a compact Riemannian manifold $M$ carries a smooth right $G$-action by isometries, with $G$ a connected Lie group, so that the Laplace–Beltrami operator commutes with the action and $L^2(M)$ decomposes into unitary representations of $G$. The associativity identity (in a real basis)
$$\sum_k M^{i,j,k}M^{\ell,m,k}=\sum_k M^{j,\ell,k}M^{i,m,k}$$
is an identity in the algebra of functions on $M$. In a $G$-adapted eigenbasis each factor sits in an isotype, and the same identity becomes a crossing equation: the two association trees for a four-linear pairing are related by the recoupling kernels of $G$,
\label{rem:xing}

\begin{equation}
\sum_{\nu}\lVert M_{\lambda\mu\to\nu}\rVert^{2}_{\mathrm{HS}}\,K(\lambda,\mu;\kappa,\rho\mid\nu)
=\sum_{\nu'}\lVert M_{\lambda\rho\to\nu'}\rVert^{2}_{\mathrm{HS}}\,K'(\lambda,\rho;\mu,\kappa\mid\nu').
\tag{\(\dagger\)}
\end{equation}

The second kernel is $K'(\lambda,\rho;\mu,\kappa\mid\nu')$, not $K'(\lambda,\kappa;\mu,\rho\mid\nu')$. With $\lambda=j_1$, $\mu=j_2$, $\kappa=j_4$, and $\rho=j_3$, this is the $K'$ of (24). The weight on the right is the multiplication map for the pair $(\lambda\rho)$, which is the pair actually crossed.
In particular, the indices in $(\dagger)$ run over irreducible representations, together with multiplicity labels when an isotype occurs more than once. They are not the indices of a chosen basis. Each Hilbert–Schmidt weight and each reduced kernel is independent of that choice.

If $G$ is compact, $K$ and $K'$ are finite recoupling coefficients, and $(\dagger)$ is an algebraic identity among the reduced elements of $\Phi$. For $SU(2)$ they are $6j$ symbols when they are three-spin associators, or the $J=0$ collapse, and $9j$ symbols when four external spins couple to a total spin $J\neq 0$.

If $G=\mathrm{PSL}_2(\mathbb{R})$ or $\mathrm{SO}_0(1,d)$ acting on $\Gamma\backslash G$, the labels run over principal series, $K,K'$ are the Plancherel recoupling kernels of Kravchuk–Mazáč–Pal, and $(\dagger)$ is the conformal bootstrap crossing system. A linear functional $\alpha$ with $\alpha\cdot K\ge 0$ on the support of the spherical Plancherel measure then yields
$$\sum_{\nu\in\Omega}\lVert M_{\lambda\mu\to\nu}\rVert_{\mathrm{HS}}^2\le C(\lambda,\mu,\Omega).$$
In particular the support of $\Phi$ in the Weyl chamber cannot place too much mass outside a compact region determined by the external labels. Those inequalities are the spectral-gap and triple-overlap bounds of Bonifacio, Kravchuk–Mazáč–Pal, and subsequent work \cite{Bon22} \cite{Bon22b} \cite{KMP24} \cite{BMP25}. Adve’s converse \cite{AA25} states that a discrete solution of $(\dagger)$ for $G=\mathrm{PSL}_2(\mathbb{R})$ arises from the unit tangent bundle of a compact hyperbolic 2-orbifold.

Crossing is associativity written in a $G$-isotypic basis. It constrains the spectrum; it does not reconstruct the metric from two given manifolds, and it does not replace matching the full table $M^{i,j,k}$ in Theorem 1.
\end{remark}

\begin{remark}[Associativity of Wigner $D$-functions and $6j$ symbols]
On $SU(2)$ take the $L^2$-normalized Peter--Weyl basis $\sqrt{2j+1}\,D^j_{mn}$ relative to the Euler measure $dg=\sin\beta\,d\alpha\,d\beta\,d\gamma$ of total mass $8\pi^2$. The structure constants are
\begin{equation}
\begin{aligned}
&M^{(j_1,m_1,n_1),(j_2,m_2,n_2),(j_3,m_3,n_3)}\\
&\qquad=\sqrt{\frac{(2j_1+1)(2j_2+1)(2j_3+1)}{8\pi^2}}\,(-1)^{m_3-n_3}
\begin{pmatrix} j_1 & j_2 & j_3 \\ m_1 & m_2 & -m_3 \end{pmatrix}
\begin{pmatrix} j_1 & j_2 & j_3 \\ n_1 & n_2 & -n_3 \end{pmatrix}.
\end{aligned}
\label{eq:su2-triple}
\end{equation}
The third columns carry minus signs because $M$ conjugates the third factor and $\overline{D^j_{mn}}=(-1)^{m-n}D^j_{-m,-n}$.
Write $I=(j,m,n)$ for a Peter--Weyl index. Associativity of pointwise multiplication in this complex basis is
\begin{equation}
\sum_K M^{IJK}\,M^{KLP}=\sum_{K'} M^{JLK'}\,M^{IK'P}.
\label{eq:assoc-complex}
\end{equation}
The contracted index sits in the first slot of the second factor. The real-basis form with that index in the third slot is not the $SU(2)$ identity.
Each $M$ factors into a left-magnetic $3j$ symbol and a right-magnetic $3j$ symbol. Summing the intermediate magnetic numbers of $K=(j,m,n)$ therefore splits as a product of two recouplings. The left-magnetic sum is Racah's identity. The sum runs over the intermediate magnetic number $m$, and $m'=-m_1-m_3$ is fixed by the selection rule:
\begin{equation}
\begin{aligned}
&\sum_m(-1)^{j-m}
\begin{pmatrix} j_1 & j_2 & j \\ m_1 & m_2 & m \end{pmatrix}
\begin{pmatrix} j & j_3 & j_4 \\ -m & m_3 & m_4 \end{pmatrix}\\
&\quad=\sum_{j'}(-1)^{j_1+j_2+j_3+j_4}\sqrt{(2j+1)(2j'+1)}
\begin{Bmatrix} j_1 & j_2 & j \\ j_4 & j_3 & j' \end{Bmatrix}\\
&\qquad\times
\begin{pmatrix} j_1 & j_3 & j' \\ m_1 & m_3 & m' \end{pmatrix}
\begin{pmatrix} j' & j_2 & j_4 \\ -m' & m_2 & m_4 \end{pmatrix}.
\end{aligned}
\label{eq:racah-left}
\end{equation}
The right-magnetic sum is the same identity with $(n_i,n,n')$ in place of $(m_i,m,m')$, and $n'=-n_1-n_3$. The product of the two recouplings still carries the external $3j$ symbols. Contracting those by $3j$ orthogonality, and summing on $j$, gives
\begin{equation}
\sum_j(2j+1)
\begin{Bmatrix} j_1 & j_2 & j \\ j_4 & j_3 & j' \end{Bmatrix}
\begin{Bmatrix} j_1 & j_2 & j \\ j_4 & j_3 & j'' \end{Bmatrix}
=\frac{\delta_{j'j''}}{2j'+1},
\label{eq:6j-orth}
\end{equation}
not an equality of two distinct $6j$ symbols, and not the Biedenharn--Elliott pentagon.

A sum over only $m_5,m_6$, with $m_1,\ldots,m_4$ held fixed, is a single admissible term. For the unsigned product the row sums freeze $m_5=-m_1-m_2$ and $m_6=-m_2-m_3$. For the signed contraction, with $-m_6$ in the contracted column, the second value is $m_6=m_2+m_3$. Either way it does not equal a bare $6j$ symbol.

The bare symbol is the full contraction
\begin{equation}
\begin{aligned}
&\begin{Bmatrix} j_1 & j_2 & j_5 \\ j_3 & j_4 & j_6 \end{Bmatrix}
=\sum_{m_1,\ldots,m_6}(-1)^{\sum_{k=1}^{6}(j_k-m_k)}
\begin{pmatrix} j_1 & j_2 & j_5 \\ m_1 & m_2 & m_5 \end{pmatrix}\\
&\qquad\times
\begin{pmatrix} j_1 & j_4 & j_6 \\ m_1 & m_4 & -m_6 \end{pmatrix}
\begin{pmatrix} j_3 & j_2 & j_6 \\ -m_3 & m_2 & m_6 \end{pmatrix}
\begin{pmatrix} j_3 & j_4 & j_5 \\ m_3 & -m_4 & m_5 \end{pmatrix}.
\end{aligned}
\label{eq:6j-def}
\end{equation}
The kernels in $(\dagger)$ are three-spin associators. An odd permutation of columns multiplies by $(-1)^{j_1+j_2+j_3+j_4+j+J}$. With the column order of \eqref{eq:racah-left},
\begin{equation}
\begin{aligned}
K(j_1,j_2;j_4,j_3\mid j)
&=(-1)^{j_1+j_2+j_3+j_4}\sqrt{(2j+1)(2J+1)}
\begin{Bmatrix} j_1 & j_2 & j \\ j_4 & j_3 & J \end{Bmatrix},\\
K'(j_1,j_3;j_2,j_4\mid j')
&=(-1)^{j_1+j_2+j_3+j_4}\sqrt{(2j'+1)(2J+1)}
\begin{Bmatrix} j_1 & j_3 & j' \\ j_4 & j_2 & J \end{Bmatrix}.
\end{aligned}
\label{eq:kernels}
\end{equation}
They vanish unless all four triangles of each symbol hold,
\begin{equation}
\begin{aligned}
&(j_1,j_2,j),\ (j,j_3,j_4),\ (j_1,j_4,J),\ (j_2,j_3,J),\\
&(j_1,j_3,j'),\ (j',j_2,j_4),\ (j_1,j_4,J),\ (j_2,j_3,J),
\end{aligned}
\end{equation}
and each triad sums to an integer. A Weyl move from the pairing $(12)(34)$ to $(13)(24)$ sends $K$ to $K'$ up to $(-1)^{j_1+j_2+j_3+j_4+j+J}$.
If $j_1,j_2,j_3,j_4$ are four external spins and $J$ is their total spin, the crossing kernel is a $9j$ symbol. The case $J\neq 0$ is the $9j$ case. A single $6j$ symbol is the associator \eqref{eq:kernels}, or the singlet collapse below.
Summing the magnetic indices in \eqref{eq:su2-triple} gives
\begin{equation}
\|M_{j_1 j_2\to j}\|^2_{\mathrm{HS}}
=\frac{(2j_1+1)(2j_2+1)(2j+1)}{8\pi^2}
\label{eq:hs}
\end{equation}
when $(j_1,j_2,j)$ is admissible, and $0$ otherwise. Substituting (24) and (26) into $(\dagger)$ gives the crossing identity
$$\begin{aligned}
&\sum_j(2j_2+1)(2j+1)^{3/2}
\begin{Bmatrix} j_1 & j_2 & j \\ j_4 & j_3 & J \end{Bmatrix}\\
&\qquad=\sum_{j'}(2j_3+1)(2j'+1)^{3/2}
\begin{Bmatrix} j_1 & j_3 & j' \\ j_4 & j_2 & J \end{Bmatrix},
\end{aligned}$$
after cancelling the common factor $(2j_1+1)(-1)^{j_1+j_2+j_3+j_4}\sqrt{2J+1}/(8\pi^2)$.

For the four-point singlet $J=0$, in the bottom-row order of the printed kernel,
\begin{equation}
\begin{Bmatrix} j_1 & j_2 & j \\ j_4 & j_3 & 0 \end{Bmatrix}
=\delta_{j_1 j_3}\,\delta_{j_2 j_4}
\frac{(-1)^{j_1+j_2+j}}{\sqrt{(2j_1+1)(2j_2+1)}},
\label{eq:6j-zero}
\end{equation}
so
\begin{equation}
K(j_1,j_2;j_4,j_3\mid j)
=\delta_{j_1 j_3}\,\delta_{j_2 j_4}\,
(-1)^{j_1+j_2+j}
\sqrt{\frac{2j+1}{(2j_1+1)(2j_2+1)}},
\label{eq:K-singlet}
\end{equation}
and likewise for $K'$ after the paired legs are swapped. Then $(\dagger)$ reduces to a finite equality of dimension factors.
\end{remark}

\section{Example}
\label{sec:ex}

Let $\set{\alpha_i} \subset \mathbb{R}^n$ be an indexed, rank $n$ lattice of Lie Algebra weights for the quotient space representation of $\mathfrak{g}=\mathbb{R}^n$ as translation invariant (i.e., constant) vector fields on itself, when $\mathbb{R}^n$ is also viewed as $\mathfrak{g}$'s associated Lie Group over a torus defined by $\mathbb{R}^n/A\mathbb{Z}^n, A \in GL(n,\mathbb{R})$. These weights define integrable lifts of 1-forms over the torus that integrate to linear functionals $\bra{x} \alpha_i\rangle,\, x\in\mathbb{R}^n$ as its Lie Group (covering the torus). These linear functionals can then be uniformly rescaled (by $2\pi \sqrt{-1}$) and exponentiated to form multiplicative characters that descend to form an orthonormal basis of $L^2(\mathbb{R}^n/A\mathbb{Z}^n,dx)$, with Lebesgue (Haar) measure $dx$.

Moreover, this basis simultaneously diagonalizes the flat torus's Laplacian, because the Laplacian is the image of a symmetric, negative-definite quadratic Casimir element under this (constant coefficient linear differential operator) quotient space representation of the universal enveloping algebra. Hence, its eigenvalues are in constant proportion (of $4\pi^2$) to the Casimir-element-determined-length-squared of each character's weight in the lattice. It is easy to see that any choice of non-degenerate, negative-definite quadratic Casimir element generates a Riemannian Geometry that is isometric to the (negative) Euclidean one.

We presently view the above basis

\begin{equation}
\{e^{2\pi\sqrt{-1}\bra{x}\alpha_i\rangle}/\sqrt{|\det A|}\}_{i=0}^\infty
\end{equation}
to be our Theorem-applicable Fourier basis of orthonormal (multiplicative character) eigenfunctions (of this quotient representation of the (negative) Euclidean Casimir element) directly corresponding to $\set{\alpha_i}$. By our Theorem's hypotheses, we must have $i < j \implies \lVert\alpha_i\rVert\leq \lVert\alpha_j\rVert$ (with the Euclidean norm on the weights).

Now we can compute

\begin{equation}
\label{eqn:a}
M^{i,j,k} = \begin{cases}
1/\sqrt{|\det A|} & \alpha_i + \alpha_j - \alpha_k = 0 \\
0 & \text{otherwise}
\end{cases}
\end{equation}
As this Equation~\ref{eqn:a} only depends on the weight lattice itself, it is orthonormal-basis-index invariant. Further, it is only invariant under linear transformations on the weight lattice $(A^{-1})^t\mathbb{Z}^n = \set{\alpha_i}$, so only an $L^2$ orthonormal eigenfunction basis map \textit{which is induced from a volume-preserving invertible linear map between two such indexed, rank $n$ weight lattices} will keep the "algebraic/topological" indexed data set $\set{M^{i,j,k}}$ invariant.

However, in order to apply Theorem~\ref{thm}, it is essential that such a linear map $B$ be $B\in SO(n,\mathbb{R})$ on the weight lattice, because the induced $L^2$ eigenfunction basis map

\begin{equation}
\{e^{2\pi\sqrt{-1}\bra{x}B\alpha_i\rangle}/\sqrt{|\det A|}\}_{i=0}^\infty
\end{equation}
must also preserve the "analytic" invariants -- the Casimir-element induced figure $\lambda_i = 4\pi^2\lVert\alpha_i\rVert^2$ for each indexed weight, i.e.\ the individual eigenvalues of the flat-tori's Laplacian.

As Milnor's duet exemplifies, having a map which preserves the lengths of the lattice weights is not sufficient to deduce the map is in $SO(n,\mathbb{R})$; we must also know that the map preserves weight angles. But this is a consequence of the formulae developed in Equation~\ref{eqn-15}:

\begin{equation}
-4\pi^2\bra{\alpha_i}\ket{\alpha_j}e^ie^j = de^i\cdot de^j = 2\pi^2\sum_k (\lVert\alpha_i\rVert^2 + \lVert\alpha_j\rVert^2 - \lVert\alpha_k\rVert^2)M^{i,j,k}e^k = 2\pi^2(\lVert\alpha_i\rVert^2 + \lVert\alpha_j\rVert^2 - \lVert\alpha_i + \alpha_j\rVert^2)e^ie^j
\end{equation}
The neat thing about this analysis is that we've proven there is no \textbf{linear} map between lattices that preserves the eigenvalues without the map being induced by a Riemannian isometry on the tori -- as a consequence of the Theorem, not because the explicit computations involved are simple polarization identities.

This representation-theoretical account \cite{AK01} is exactly equivalent to the prior development of \textit{lattice congruence} \cite{NRR22} traditionally used to delineate isometry classes of flat tori. In fact, the matrix transpose of such a linear map on the weight lattice $B\in SO(n,\mathbb{R})$, \textbf{is} the contravariant Riemannian isometry between the tori, as provided by application of the \textit{Gelfand-Naimark Representation Theorem} during the proof of Theorem~\ref{thm}.

\paragraph{Acknowledgements}
The original research was funded in part by a gracious James Simons Research Award in 1995-1996, and the generous support of an Alfred P. Sloan Dissertation Fellowship in 1996-1997 at the University at Stony Brook.

The author would also like to thank Tanya Christiansen, Carolyn Gordon, Hamid Hezari, Harish Seshadri, and especially Leon Takhtajan for their technical assistance and review in the preparation of this manuscript for publication.

\doclicenseThis

\printbibliography

\end{document}